\documentclass[10pt,a4paper,twoside]{amsart}

\usepackage[all]{xy}
\xyoption{poly}
\CompileMatrices
\usepackage{amsmath,amssymb,latexsym}
\usepackage{graphics,pstricks,textcomp}
\usepackage{mathrsfs}

\def\bc{\begin{center}}
\def\ec{\end{center}}
\def\be{\begin{equation*}}
\def\ee{\end{equation*}}
\def\ba{\begin{array}}
\def\ea{\end{array}}
\def\bp{\begin{picture}}
\def\ep{\end{picture}}

\def\bea{\begin{eqnarray*}}
\def\eea{\end{eqnarray*}}
\def\bnea{\begin{eqnarray}}
\def\enea{\end{eqnarray}}

\def\lra{\longrightarrow}

\def\mt{\mapsto}
\def\lmt{\longmapsto}
\def\bc{\begin{center}}
\def\ec{\end{center}}
\def\noi{\noindent}

\def\A{\mathbf{A}}

\def\C{\mathbb{C}}
\def\N{\mathbb{N}}
\def\R{\mathbb{R}}
\def\Z{\mathbb{Z}}

\def\P{\mathscr{P}}

\def\ob{\mathrm{Ob}\,}
\def\ad{\mathrm{ad}}

\def\sl{{\mathfrak{sl}}_2}
\def\RG{\mathbf{Rep}_f (\sl)}
\def\id{\mathrm{Id}}

\def\de{{\bf Definition} : }

\def\lem{{\bf Lemma} : }
\def\pr{{\it Proof} : }
\def\ie{{\it i.e. }\,}

\def\TT{\widetilde{T}}

\def\FT{\widetilde{F}}

\def\LT{\widetilde{L}}

\def\f{\varphi}
\def\F{\Phi}
\def\g{\gamma}
\def\G{\Gamma}
\def\s{\sigma}
\def\sb{\overline{\sigma}}

\def\d{\delta}
\def\e{\varepsilon}
\def\p{\partial}

\def\tr{\mathrm{tr}}

\def\homa{\hom_{\A}}

\def\pic{\xymatrix}
\def\vs{\vspace{5mm}}
\def\nl{\newline}

\begin{document}

\fontfamily{cmss}
\fontseries{m}
\selectfont

\hfill\today

\title{Combinatorial Stacks and the Four-Colour Theorem}

\author{Romain Attal}

\email{attal@lpthe.jussieu.fr}

\address{Laboratoire de Physique Th\'{e}orique et Hautes \'{E}nergies \nl
Universit\'e Pierre et Marie Curie \nl
4, Place Jussieu \nl
F-75005 Paris (FRANCE)}

\begin{abstract} 
We interpret the number of good four-colourings of the faces of a trivalent, spherical polyhedron as the 
2-holonomy of the 2-connection of a fibered category, $\f$, modeled on $\RG$ and defined over the 
dual triangulation, $T$. We also build an $\sl$-bundle with connection over $T$, that is a global, 
equivariant section of $\f$, and we prove that the four-colour theorem is equivalent to the fact 
that the connection of this $\sl$-bundle vanishes nowhere. This interpretation is proposed as a first 
step toward a cohomological proof of the four-colour theorem.
\end{abstract}

\maketitle

\tableofcontents


\section{Introduction}

Let us consider a finite spherical polyhedron, $P$, and a palette of four colours, $\{W,R,G,B\}$.
We will call a good colouring of $P$ any map which associates one of these colours to each face of $P$ in such a
way that any two adjacent faces carry distinct colours. The four-colour theorem \cite{AC,AH,RSST} states that
such a map exists for any $P$. The goal of the present work is to provide a geometric interpretation of this 
theorem. We obtain here two new results : the number of good colourings of a trivalent, spherical 
polyhedron is the 2-holonomy of a 2-connection on a fibered category over the dual triangulation, 
$T=P^\ast$ (Theorem 2) ; the four-colour theorem is equivalent to the existence of a non-vanishing, 
equivariant global section of this fibered category (Theorem 4). 

In order to study the colourability of $P$, let us start by making some classical modifications.
We first remark that it is sufficient to prove the colourability of trivalent polyhedra. 
Indeed, by cutting a little disk around each vertex of degree $>3$ in $P$, one obtains 
a trivalent polyhedron and each good colouring of the latter provides a good 
colouring of $P$ by shrinking this disk to the initial vertex. Henceworth, we will 
suppose that $P$ is itself trivalent. Secondly, let us identify our four colours 
with the pairs of diametrally opposite vertices of a cube : 
$W=\{w,w'\}$, $R=\{r,r'\}$, $G=\{g,g'\}$ and $B=\{b,b'\}$. Then each good colouring of the 
three faces which surround a vertex of $P$ defines an edge-loop in this cube such that the 
determinant of any triplet of successive vectors be $\pm 1$ :

\subsubsection*{}
\subsubsection*{}
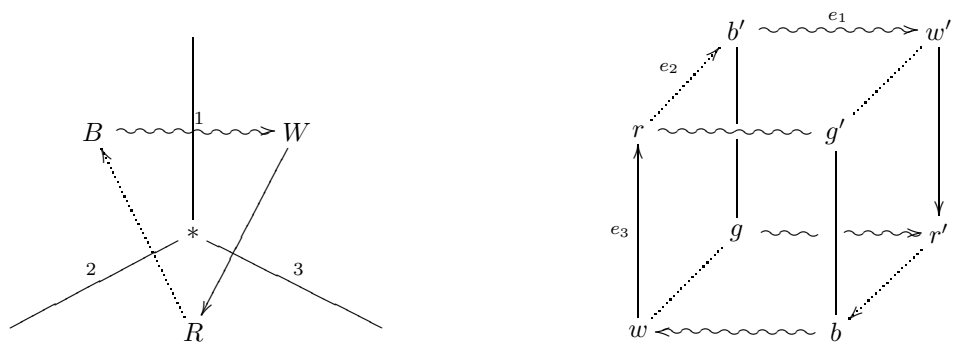
$$
\pic{
& & \ar@{-}[dd] & & & & & & b' \ar@{-}[d] \ar@{~>}[rr]^{e_1} & & w' \ar@{.}[dl] \ar[dd] \\
& B \ar@{~>}[rr]^{\ 1} & & W \ar[ddl] & & & & r \ar@{.>}[ur]^{e_2} \ar@{~}[rr] & \ar@{-}[d] & 
g' \ar@{-}[dd] & \\
& & \ast \ar@{-}[dll]_{2} \ar@{-}[drr]^{3} & & & & & & 
g \ar@{~}[r] & \ar@{~>}[r] & r' \ar@{.>}[dl] \\
& & R \ar@{.>}[uul] & & & & & w \ar@{->}[uu]^{e_3} \ar@{.}[ur] \ar@{<~}[rr] & & b & \\
}
$$
\subsubsection*{}
\subsubsection*{}

\noi A map $(u:T_1\to\{e_1,e_2,e_3\})$ satisfying this property will be called a good numbering.
Thus, the number of good numberings of the edges of $T$ is one quarter of the number of good 
colourings of the faces of $P$, as proved by P.G. Tait \cite{T}. We call this integer, $K_T$, the chromatic 
index of $T$ and the four-colour theorem states that $K_T\neq 0$ for any finite, spherical triangulation, $T$.

Our article is organised as follows.
In Section 2, we give a proof of Penrose's formula which expresses $K_T$ as a partition function.
In Section 3, we define the graph $\P$ of edge-paths of $T$.
In Section 4, we collect useful results about representations of $\sl$.
In Section 5, we construct the chromatic stack, $\f$, which is a fibered category over $T$, endowed 
with a functorial 1-connection and with a natural 2-connection, and we prove that $K_T$ is the 
2-holonomy of this 2-connection on $T$.
In Section 6, we define another fibered category, $\Phi$, over $\P$. 
By integrating the functorial connection of $\Phi$ along a 2-path which sweeps each triangle of $T$ 
once only, we obtain an equivariant global section of the pull-back of the chromatic stack to a triangulation 
$\TT$ of the disk. This section, $\zeta$, is an $\sl$-bundle with connection whose holonomy on $\p\TT$ is $K_T$.
Our construction is an adaptation of Stokes theorem to a case of combinatorial differential forms with values 
in the tensor category $\A=\RG$ and we can write it symbolically $K_T=\int_T\f=\int_{\p\TT}\zeta$.
Since $K_T$ depends linearly on the value of $\zeta$ on each inner edge of $\TT$, we obtain this way 
our second result : the four-colour theorem is equivalent to the fact that $\zeta$ vanishes nowhere.

\section{The chromatic index}

The idea to translate the four-colour problem in terms of linear algebra 
is due to Roger Penrose. Let us fix a finite, spherical triangulation $T=(T_0,T_1,T_2)$. 
$T_0$ is the set of its vertices, $T_1$ the set of its edges and $T_2$ the set of its triangles. 
Following \cite{RP}, we define the chromatic index of $T$ as

\subsubsection*{}
\bnea
\boxed{K_T := \sum_u \prod_{[xyz]} i\,\det\,(u_{xy},u_{yz},u_{zx})} 
\enea
\subsubsection*{}
\subsubsection*{}

\noi In this sum, $u$ runs over the set of all maps from $T_1$ to $\{e_1,e_2,e_3\}$, 
the canonical basis of $\R^3$, and $[xyz]$ runs over the set of positively oriented 
triangles of $T$. The integrality of $K_T$ follows from the fact that, if $u$ is a good numbering of $T_1$, 
\ie if no determinant vanishes in this product, the number of triangles where $\det=(+1)$ minus 
the number of triangles where $\det=(-1)$ is a multiple of $4$, as proves the following lemma.

\subsubsection*{}
\lem \emph{If $u$ is a good numbering of $T_1$ and if $n_+$ 
(resp. $n_-$) denotes the number of triangles $[xyz]$ such that 
$\det (u_{xy},u_{yz},u_{zx}) = (+1)$ (resp. $(-1)$), then $n_+\equiv n_-$ mod 4.}
\subsubsection*{}

\pr 1) Starting from $(T,u)$, we can build another triangulation, $T'$, equipped with a good 
edge numbering, $u'$, such that ${n'}_+=0$. Indeed, if two adjacent, positively oriented 
triangles of $T$, say $[xyz]$ and $[zyw]$, have $\det=(+1)$ (positive triangles), then we 
can flip their common edge $[yz]$ to $[xw]$ and obtain a new pair of negative triangles, 
$[xyw]$ and $[wzx]$, where $\det = (-1)$. During this step, $(n_+-n_-)$ is reduced by $4$. 
Once all these pairs of neighbour positive triangles have been eliminated this way, 
the remaining contributions to $n_+$ are triangles surrounded by three negative neighbours. 
By adjoining three edges and a trivalent vertex inside each isolated triangle of this kind, 
we change a positive triangle for three negative ones. Again, $(n_+-n_-)$ is reduced by $4$, 
and $(T',u')$ is reached at the end of this process.

\subsubsection*{}
2) Consider all pairs of triangles, $[xyz]$ and $[zyw]$, with ${u'}_{yz}=e_1$ on their 
common edge, $[yz]$. Since $\det (u_{xy},u_{yz},u_{zx}) = \det (u_{zy},u_{yw},u_{wz}) = (-1)$, 
the opposite sides of the rectangle $[xywz]$ carry the same vector, say ${u'}_{xz}={u'}_{yw}=e_2$ 
and ${u'}_{xy}={u'}_{zw}=e_3$. Let us join the midpoints of two opposite edges with a simple curve. 
By repeating this process inside all such pairs of triangles, we obtain two simple closed curves, $c_2$ 
and $c_3$. If we orient these curves conveniently, their intersection number is equal to 
$\vert{u'}^{(-1)}(e_1)\vert$, the number of edges of $T'$ marked with $e_1$. But, after Jordan's theorem, 
the intersection number of two simple closed curves in $S^2$ is even. Therefore, $\vert{u'}^{(-1)}(e_1)\vert$,
the number of edges mapped to $e_1$ by $u'$, is even. Similarly, $\vert{u'}^{(-1)}(e_2)\vert$ 
and $\vert{u'}^{(-1)}(e_3)\vert$ are also even, as well as the total number of edges of $T'$ :

\bea
{t'}_1=\vert{u'}^{(-1)}(e_1)\vert+\vert{u'}^{(-1)}(e_2)\vert+\vert{u'}^{(-1)}(e_3)\vert \in 2\,\N \\
\eea

3) Since $T'$ is a triangulation of a closed surface, we have $3\,{t'}_2=2\,{t'}_1$. Since 
${t'}_1$ is even, we obtain ${t'}_2 = {n'}_- \in 4\,\N$. Therefore, $n_+$ and $n_-$ are 
congruent modulo 4 :

\bnea
\boxed{ (n_+-n_-) \in 4\,\Z } 
\enea
\hfill\qed

\subsubsection*{}
${\bf Theorem \  1}$ : \emph{$K_T$ is the number of good numberings of $T_1$.} 
\subsubsection*{}

\pr If $u$ is a bad numbering, then one of the determinants is zero and the corresponding 
product vanishes. On the other hand, if $u$ is a good numbering, then the corresponding 
product is equal to $i^{(n_+-n_-)} = 1$, after the precedent lemma. Therefore, the sum of all 
these products equals the number of good numberings of $T_1$.

\hfill\qed

\subsubsection*{}

\section{The graph of edge-paths}

Having fixed our triangulation, $T$, let us define the graph $\P$ whose vertices 
are the edge-paths of $T$ :

\bea
\P_0 = \bigcup_{\ell\geqslant 0} \big\lbrace \g=(x_0,\cdots,x_\ell) \ :
\ \{x_i,x_{i+1}\} \in T_1 \ \forall\, i \big\rbrace \\
\eea
and whose edges, called the 2-edges of $T$, are the pairs of paths, with 
the same source and the same target, which bound a single triangle of $T$ :

\bea
\P_1 = \big\lbrace \{ (x_0,\cdots,x_\ell),(x_0,\cdots,x_i,y,x_{i+1},\cdots,x_\ell) \}
\in \P_0\times\P_0 \ : \   \{x_i y x_{i+1}\} \in T_2 \big\rbrace \\
\eea

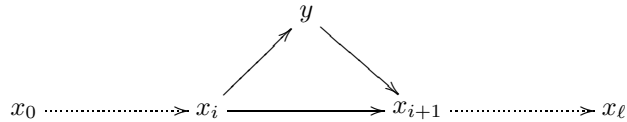
$$
\pic{
& & & y \ar@{->}[dr] & & & \\
x_0 \ar@{.>}[rr] & & x_i \ar@{->}[rr] 
\ar@{->}[ur] & & x_{i+1} \ar@{.>}[rr] & & x_\ell \\
& & & & & & \\
}
$$

The oriented 2-edges are the corresponding ordered pairs.
A 2-path in $T$ is an edge-path in $\P$, \ie a family $\G=(\g_0,\cdots,\g_n)$ 
such that $\{\g_i,\g_{i+1}\}\in\P_1$ for $i=0,\cdots,n-1$. They form the set $\P_2$ :

\bea
\P_2 = \bigcup_{n\geqslant 0} \big\lbrace \G=(\g_0,\cdots,\g_n) \ : 
\ \{\g_i,\g_{i+1}\}\in\P_1 \quad \text{for} \quad i=0,\cdots,n-1 \\
\eea
For each 2-path $\G=(\g_0,\cdots,\g_n)$, there is a 2-path $\widetilde{\G}$ 
going backward in time :

\bea
\widetilde{\G}=(\g_n,\cdots,\g_0) \\
\eea
The 0-source (resp. 0-target) of $\G$ is the common source (resp. target) of 
the $\g_i$'s. The 1-source of $\G$ si $\g_0$ and and its 1-target is $\g_n$.
The oriented 2-cells of $T$ are its smallest 2-paths. 
They have the form $\big( (xz),(xyz) \big)$ or $\big( (xyz),(xz) \big)$, 
for some triangle $\{xyz\}$.

\section{Representations of $\sl$}

As we have seen above, Penrose's formula involves the determinants of triples of basis 
vectors of $\R^3$. If we endow $\R^3$ with its canonical euclidian structure and with the 
corresponding cross-product, we obtain a Lie algebra isomorphic to ${\mathfrak{so}}_3$. 
Since we will use complex coefficients and Schur's lemma, valid only for representations over an 
algebraically closed field, we will work with its complexification, $V=\sl$. 
We will note $I=\id_{V}$, $V^\ell=V^{\otimes\ell}$  and  $I^\ell=\id_{{V}^{\ell}}$, 
where $V^\ell$ carries the representation

\bea
\rho_\ell &:& V \lra {\mathrm{End}} ({V}^{\ell}) \\
& & x \lmt \rho_\ell (x) = \sum_{k=1}^\ell I^{k-1} \otimes \ad_x \otimes I^{\ell-k} \\
\eea
Let $\A=\RG$, the category of finite dimensional representations of $\sl$ over complex 
vector spaces. 
If $M$ and $M'$ are two $V$-modules, carrying, respectively, the representations $R$ and $R'$, 
we will often identify $M$ with $M\otimes -$, the endofunctor of $\A$, and write $M'M$ for $M'\otimes M$.
For each $j\in\frac 12 \N$, let $(R_j:V \to {\mathrm{End}} (V_j))$ be a representative 
of the isomorphy class of representations of spin $j$ and dimension $2j+1$. For example, 
we can choose $V_0=\C$, $V_{1/2}=\C^2$ and $V_1=V$. After Schur's lemma, the irreducible 
representations are orthonormal for the bilinear bifunctor $\homa$ :

\bnea
\boxed{\homa (R_j,R_k) \simeq \delta_{jk} R_0} 
\enea
\subsubsection*{}
\noi The intertwining number between two representations $R$ and $R'$ is defined as the 
dimension of the space $\homa (R,R')$ :

\bea
c(R,R') = \dim_{\C} \big( \homa (R,R') \big) \\
\eea
After Clebsch-Gordan's rule, $V^2\simeq V_0\oplus V_1\oplus V_2$ and $c(V,V^2)=c(V^2,V)=1$.
The projectors onto the isotypic components of $V^2$, of spin $0$, $1$ and $2$, respectively 
map $u\otimes v$ to

\bea
T (u\otimes v) &=& (u \cdot v) \, e_a \otimes e_a \\
A (u\otimes v) &=& \frac 12 (u_a v_b - u_b v_a) \, e_a \otimes e_b \\
S (u\otimes v) &=& \frac 12 (u_a v_b + u_b v_a) \, e_a \otimes e_b - (u_a v_a) \, e_a \otimes e_a \\
\eea
The line $L=\homa(V,V^2)$ is spanned by the map $F$ defined by

\bea
F(e_a)=i \, e_{a-1}\wedge e_{a+1} \\
\eea
and the line $\LT=\homa(V^2,V)$ is spanned by the bracket, noted $\FT$ :

\bea
\FT(e_a\otimes e_b)=[e_a,e_b]=i\,\e_{abc} \, e_c \\
\eea
All these morphisms of representations satisfy the relations

\bea
\FT F &=& 2\, I \\
F \FT &=& A \\
T+A+S &=& I^2 \\
(\FT \otimes I) (I \otimes F) & = & (I \otimes \FT) (F \otimes I) \\
&=& T+2A-2S \\
F &=& (\FT \otimes I) (I \otimes F) F  \\
\FT &=& \FT (I \otimes \FT) ( F \otimes I)  \\
\eea
\subsubsection*{}

\section{The chromatic stack, $\f$}

The notion of combinatorial stack appeared in \cite{Kap} and we used it in \cite{A} 
to give a construction of non-abelian $G$-gerbes over a simplicial complex. Dually, 
we can also use coefficients in a category of representation. Thus, we define the chromatic 
stack, $\f$, as a 2-functor which represents the simplicial homotopy groupoid $\Pi_1(\P)$ 
into the 2-category of $\A$-modules.  $\f$ is generated by pasting the following data :

\bea
\f_x &=& \A \\
\f_{xy} &=& (V \otimes - : \f_y \to \f_x ) \\
\f_{(x_0,\cdots,x_\ell)} &=& (V^\ell \otimes - : \f_{x_\ell} \to \f_{x_0}) \\
\f_{\s} &=& F \quad {\mathrm{if}} \quad \s=\big( (xyz),(xz) \big) \\
&=& \FT \quad {\mathrm{if}} \quad \s=\big( (xz),(xyz) \big) \\
\f_{\g\g'} &=& (I^k\otimes\f_{\s}\otimes I^{\ell-k-1} : \f_{\g'} \to \f_{\g}) \\
\f_{(\g_0,\cdots,\g_n)} &=& (\f_{\g_0\g_1}\circ\cdots\circ\f_{\g_{n-1}\g_n}: \f_{\g_n}\to\f_{\g_0}) \\
\eea
The 1-connection of $\f$ is the family of functors $(\f_\g)_{\g\in \P_0}$, and the 2-connection 
of $\f$ is the family of natural transformations $(\f_\G)_{\G\in \P_2}$. 
In order to compute the chromatic index, we choose a 2-loop, 
$\G=(\g_0,\cdots,\g_n)$, based at $(a,b)=\g_0=\g_n$, and sweeping each triangle of $T$ 
once only. To each path $\g_p=(a,x_{p1},\cdots,x_{p\ell_{p-1}},b)$, of length 
$\vert\g_p\vert=\ell_p$, $\f$ associates a copy of $V^{\ell_p}$. For each $p\in\{2,\cdots,n\}$, 
the loop $\g_p$ differs from $\g_{p-1}$ either by the insertion of a vertex $y\in T_0$ between 
$x_{p-1,k_p}$ and $x_{p-1,k_p+1}$ or by the deletion of $x_{p-1,k_p}$, 
where $x_{p-1,k_p-1}$ and $x_{p-1,k_p+1}$ are supposed to be adjacent.
Each such move is represented by a linear map of the form

\bea
\f_{\g_{p-1} \g_{p}} &=& F_{k_p \ell_p} 
\ =\ (I_{V^{k_p-1}} \otimes F \otimes I_{V^{\ell_p-k_p}} \ : \ V^{\ell_p}\lra V^{\ell_p+1}) 
\qquad {\mathrm{if}} \quad \ell_{p-1}=\ell_p+1 \\
&=& \FT_{k_p \ell_p} 
\ =\ (I_{V^{k_p-1}} \otimes \FT \otimes I_{V^{\ell_p-k_p-1}} \ : \ V^{\ell_p}\lra V^{\ell_p-1}) 
\qquad {\mathrm{if}} \quad \ell_{p-1}=\ell_p-1 \\
\eea

Since Penrose's formula looks like the partition function of a statistical model, it is natural to 
express $K_T$ as the trace of a product of transfer matrices which represent linear maps 
between tensor powers of $V$. This approach will give us an efficient way to compute it, 
because the bad numberings are eliminated progressively during the sweeping process.
Geometrically, the construction of the chromatic stack allows us to reinterpret $K_T$ as 
a $2$-holonomy, which is the categorical analogue of a holonomy in a fiber bundle.

\subsubsection*{}
\de \emph{The 2-holonomy of $\f$ along a 2-loop $\G=(\g_0,\g_1,\cdots,\g_{n-1},\g_0)$ 
based at $\g_0$, is the natural transformation}

\bea
\f_\G = \f_{\g_0\g_1}\circ\cdots\circ\f_{\g_{n-1}\g_0} &:& \f_{\g_0} \lra \f_{\g_0} \\
\eea
When $\g_0=(a)$, $\f_\G$ is an endomorphism of $\f_{\g_0}=\id_{\A}$ so that $\f_\G$ defines 
canonically a complex number. Moreover, after the following theorem, which illustrates the 
pasting lemma \cite{Po} in the 2-category of $\A$-modules, the trace of 
$\f_\G\in{\mathrm{End}}(\f_{\g_0})$ depends only on $T$ and not on the 2-path $\G$.

\subsubsection*{}
${\bf Theorem \  2}$ : \emph{If $\G$ is a 2-loop which sweeps each triangle of $T$ once only, 
then the trace of the 2-holonomy of $\f$ along $\G$, evaluated in the 
representation associated to the base path of $\G$, is the chromatic index of $T$ :}

\bnea
\boxed{\tr_{\f_{\g_0}}(\f_\G) = K_T} 
\enea
\subsubsection*{}

\pr Let $\G=(\g_0,\g_1,\cdots,\g_{n-1},\g_0)$ be such a 2-loop. 
Let $p\in\{0,\cdots,n-1\}$ and suppose that $\g_{p+1}$ is obtained from $\g_p$ 
by inserting $y$ between $x_j$ and $x_{j+1}$, with $x_j \neq y \neq x_{j+1} \neq x_j$ :

\bea
\g_p &=& (x_0, \cdots, x_\ell ) \\
\g_{p+1} &=& (x_0,\cdots, x_j, y, x_{j+1}, \cdots, x_\ell ) \\
\eea
Then the 2-arrow $\f_{\g_p \g_{p+1}}$ is the intertwiner

\bea
\f_{\g_p \g_{p+1}} = I_{\f_{x_0 x_1}} \otimes\cdots\otimes 
I_{\f_{ x_{j-1} x_j }} \otimes \FT \otimes I_{\f_{ x_{ j+1, j+2 }}} 
\otimes\cdots\otimes I_{\f_{ x_{\ell-1} x_\ell }} = \FT_{\ell k} \\
\eea
which is represented by the matrix $M_p$ whose entries are given by

\bea
M_{p,ab} 
&=& \d_{a_0 b_0} \cdots \d_{a_{j-1} b_{j-1}} 
\, \big( i \, \e_{a_j b_j b_{j+1}} \big) \,
\d_{a_{j+1} b_{j+2}} \cdots \d_{a_{\ell-1} b_\ell} \\
\eea
If $\g_{q+1}$ is obtained from $\g_q$ by deleting a vertex between $y_k$ and $y_{k+1}$,
then $\f_{\g_q \g_{q+1}}$ is the intertwiner going backwards

\bea
\f_{\g_q \g_{q+1}} = I_{\f_{y_0 y_1}} \otimes\cdots\otimes 
I_{\f_{ y_{k-1} y_k }} \otimes F \otimes I_{\f_{ y_{ k+1, k+2 }}} 
\otimes\cdots\otimes I_{\f_{ y_{\ell-1} y_\ell }} = F_{\ell+1,k} \\
\eea
and is represented by the matrix whose entries are

\bea
M_{q,ab} &=& \d_{a_0 b_0} \cdots \d_{a_{k-1} b_{k-1}} 
\, \big( (-i) \, \e_{a_k b_{k+1} b_k} \big) \,
\d_{a_{k+2} b_{k+1}} \cdots \d_{a_\ell b_{\ell-1}} \\
\eea
Now, let $a^p = \big(a^p_1, \cdots , a^p_{\ell_p} \big)$ be a generic multi-index for the basis 
vectors of the representation $\f_{\g_p}$, with $a^p_j \in \{1,2,3\}$ for $j=1,\cdots,\ell_p$. 
The number $\tr_{\f_{\g_0}}(\f_\G)$ is the trace of the product of these matrices :

\bea
\tr_{\f_{\g_0}}(\f_\G) &=& \sum_a 
\prod_{p=0}^{n-1} M_{p,{a^p} a^{p+1}} \\
\eea
In this sum, $a$ runs over the set of families $(a^0, \cdots, a^n)$ of multi-indices
$a^p = \big( a^p_1, \cdots , a^p_{\ell_p} \big)$ with $a^p_i \in \{ 1,2,3 \}$.
To each edge of $T$ are associated as many indices as there are paths 
$\g_p$ which contain it. Let $N_{xy}$ be the number of indices associated to $(xy)$. 
Among them, $(N_{xy}-2)$ indices are constrained by the $\d$'s to be equal. 
Similarly, the two $\e$'s associated to the two triangles which contain $(xy)$ force
the two remaining indices to take the same value. Since the $\d$'s are sandwiched 
between these two $\e$'s, these two indices are in fact equal and there is one and 
only one free index $a_{xy}$ associated to each edge $(xy)$.
The various factors of the product are equal to one except for the $\e$'s which 
can be indexed by the positively oriented triangles of $T$. 
Therefore, the precedent formula becomes

\bea
\tr_{\f_{\g_0}}(\f_\G) &=& \sum_{(xy)\in T_1} \quad \sum_{a_{xy}\in\{1,2,3\}} 
\bigg( \prod_{[xyz]\in T_2} i\, \e_{a_{xy} a_{yz} a_{zx}} \bigg) \\
&=& \sum_u \prod_{[xyz]} i\,\det (u_{xy},u_{yz},u_{zx}) \\
\eea
where $u$ describes the set of all maps from $T_1$ to $\{e_1,e_2,e_3\}$ 
and the triangles $[xyz]$ all have the same orientation.

\subsubsection*{}
\hfill \qed
\subsubsection*{}

Initially, $K_T$ is defined as a sum of $3^{t_1}$ terms and most of them vanish. 
By working in the tensor algebra, $T(V)$, the bad bumberings are eliminated during 
the sweeping process and the computation is much quicker if we use formula (4).
Moreover, this method provides explicitely all good numberings.
\subsubsection*{}

{\bf Example} : Let us apply the relation (4) to the computation of the chromatic index 
of the octahedron.

$$
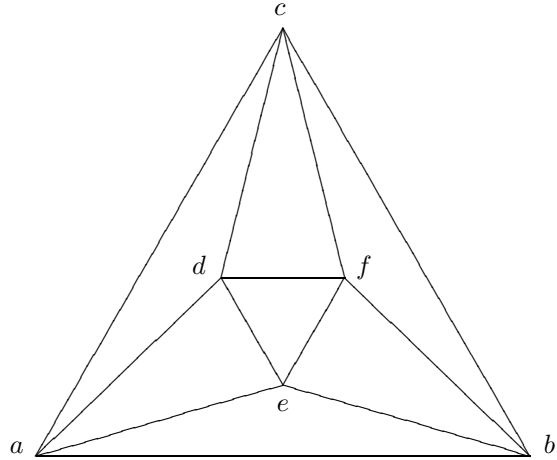
\begin{picture}(150,200)(10,-50)
\xy
/l3pc/:,{\xypolygon3"A"{~:{(.75,0):}}},
{\xypolygon3"B"{~:{(-3,0):}}},
{"A1"\PATH~={**@{-}}'"B3"'"A2"'"B1"'"A3"'"B2"'"A1"}
\endxy
\put(-197,-53){$a$}
\put(5,-53){$b$}
\put(-97,113){$c$}
\put(-128,15){$d$}
\put(-96,-37){$e$}
\put(-66,15){$f$}
\end{picture}
$$
\subsubsection*{}
\subsubsection*{}

\noi We sweep this triangulation with the $2$-path 
$\G=\big( (ab),(aeb),(adeb),(adefb),(adfb),(adcfb),(acfb),(acb),(ab) \big)$.
For simplicity, we will write ${\bf a}_1 \cdots {\bf a}_\ell$ for 
$e_{a_1}\otimes\cdots\otimes e_{a_\ell}$ with ${\bf a}_i\in\{{\bf 1,2,3}\}$. 
The successive images of ${\bf 1}$ via the maps $\f_{\g\g'}$ are :

\bea
{\bf 1} & \mt & i ({\bf 23 - 32}) \\
& \mt & i^2 ({\bf 313 - 133 -122 + 212}) \\
& \mt & i^3 ({\bf 3112 - 3121 - 1312 + 1321 - 1231 + 1213 + 2131 - 2113}) \\
& \mt & i^4 ({\bf - 331 - 122 - 111 - 111 - 133 - 221}) \\
& \mt & i^5 ({\bf - 3121 + 3211 - 1312 + 1132 - 1231 + 1321 - 1231 + 1321 - 1123 + 1213 - 2311 + 2131}) \\
& \mt & i^6 ({\bf - 221 - 111 + 212 - 331 - 221 - 331 - 221 + 313 - 111 - 331}) \\
& \mt & i^7 ({\bf 23 + 23 - 32 + 23 - 32 + 23 - 32 - 32}) \\
& \mt & i^8 ({\bf 1 + 1 + 1 + 1}) \ = \ 4 \cdot {\bf 1} \\
\eea
Consequently, $K_{octa.}=3! \cdot 4=24$ and there exist $4\cdot 24=96$ good colourings 
of the dual cube. We have made $64$ operations instead of $3^{12}=531441$. It would be 
interesting to evaluate the complexity of this method for generic triangulations.
Using the same method, one can compute the chromatic index of the icosahedron and one finds
$K_{ico.}=60$, proving this way that there exist 240 good colourings of the faces of the dual 
dodecahedron.

\section{A global section of $\f$}

$\f$ induces over $\P$ another fibered category, $\F$, defined as follows. 
To each path $\alpha=(a_0, \cdots , a_\ell)$, we associate the category $\F_\alpha$ 
whose objects are the sections of $\f$ over $\alpha$. These are the families of 
$V$-modules, $\zeta_{a_i} \in\ob(\f_{a_i})$, connected by intertwiners :

\subsubsection*{}
$$
\zeta = \bigg( 
\pic{
\zeta_{a_{i-1}} \ar@/_1pc/[rr]_{\zeta_{a_i a_{i-1}}} & & 
V\zeta_{a_i} \ar@/_1pc/[ll]_{\zeta_{a_{i-1} a_i}} \\
}
\bigg)_{1 \leqslant i \leqslant \ell} \in \ob (\F_{\alpha})
$$
\subsubsection*{}

\noi If $\zeta,\omega\in\ob(\F_\alpha)$, then $\hom_{\F_\alpha}(\zeta,\omega)$ is the 
vector space of families $(u_i:\zeta_{a_i}\to {\omega}_{a_i})_{0\leqslant i\leqslant\ell}$ 
of linear maps such that the following diagrams commute :

\subsubsection*{}
\[
\pic{
\zeta_{a_{i-1}} \ar[dd]_{u_{i-1}} \ar[rr]^{\zeta_{a_i a_{i-1}}} & & 
V \zeta_{a_i}\ar[dd]^{I \otimes u_i} & & & & 
\zeta_{a_{i-1}} \ar[dd]_{u_{i-1}} \ar@{<-}[rr]^{\zeta_{a_{i-1} a_i}} & & 
V \zeta_{a_i}\ar[dd]^{I \otimes u_i} \\
& & & & & & & & \\
{\omega}_{a_{i-1}} \ar[rr]_{{\omega}_{a_i a_{i-1}}} & & V {\omega}_{a_i} & & & & 
{\omega}_{a_{i-1}} \ar@{<-}[rr]_{{\omega}_{a_{i-1} a_i}} & & V {\omega}_{a_i} \\
}
\]

\bea
{\omega}_{a_i a_{i-1}} \circ u_{i-1} = (I \otimes u_i) \circ \zeta_{a_i a_{i-1}} 
& \hspace{25mm} & u_{i-1} \circ \zeta_{a_{i-1} a_i} = {\omega}_{a_{i-1} a_i} \circ (I \otimes u_i) \\
\eea

\de  \emph{Let $\alpha = (a_0, \cdots, a_\ell)$ be a path of length $\ell$ 
and let $\zeta\in\ob (\F_\alpha)$ be a section of $\f$ over $\alpha$. 
The direct transport operator of $\zeta$ along $\alpha$ is the morphism}

\bea
T_{\alpha}(\zeta) = (I^{\ell-1}\otimes \zeta_{a_{\ell} a_{\ell-1}}) 
\circ (I^{\ell-2}\otimes \zeta_{a_{\ell-1} a_{\ell-2}}) \circ\cdots\circ 
(I\otimes \zeta_{a_2 a_1})\circ \zeta_{a_1 a_0} \ : \ \zeta_{a_0} \lra V^{\ell} \zeta_{a_\ell} \\
\eea
\emph{and the inverse transport operator of $\zeta$ is the morphism}

\bea
\overline{T}_{\alpha}(\zeta) = \zeta_{a_0 a_1} \circ (I\otimes \zeta_{a_1 a_2})
\circ\cdots\circ (I^{\ell-2}\otimes \zeta_{a_{\ell-2} a_{\ell-1}}) \circ 
(I^{\ell-1}\otimes \zeta_{a_{\ell-1} a_{\ell}})  \ : \ 
V^\ell \zeta_{a_\ell} \lra \zeta_{a_0} \\
\eea
Let us note that $T_{\overline{\alpha}}\neq{\overline{T}}_\alpha$.
$\F_{\vert\P_1}$ is generated by its restriction to the oriented 2-cells of $T$. 
If $\s=\big( (xz),(xyz) \big)$ and if $\zeta\in\ob(\F_{(xyz)})$, then we define 
${\xi}=\F_{\s}(\zeta)\in\ob(\F_{(xz)})$ by

\bea
{\xi}_x &=& {\zeta_x} \\
{\xi}_z &=& {\zeta_z} \\
{\xi}_{zx} &=& (\FT \otimes I_{\zeta_z}) \circ (I \otimes {\zeta_{zy}}) \circ {\zeta_{yx}} \\ 
{\xi}_{xz} &=& {\zeta_{xy}} \circ (I \otimes \zeta_{yz}) \circ (F \otimes I_{\zeta_z}) \\ 
\eea
Similarly, if ${\xi}\in\ob(\F_{(xz)})$, we define $\zeta=\F_{\sb}({\xi})\in\ob(\F_{(xyz)})$ by

\bea
\zeta_x &=& {\xi}_x \\
\zeta_y &=& V {\xi}_z \\
\zeta_z &=& {\xi}_z \\
\zeta_{yx} &=& (F \otimes I_{\xi_z}) \circ \xi_{zx} \\
\zeta_{xy} &=& \xi_{xz} \circ (\FT \otimes I_{\xi_z}) \\
\zeta_{yz} &=& I\otimes I_{\xi_z} \ = \ \zeta_{zy} \\
\eea
If $u\in\hom_{\F_{(xyz)}} (\zeta,\omega)$, then we have the commutative diagrams

\subsubsection*{}
$$
\pic{
\zeta_x \ar[dd]_{u_x} \ar[rr]^{\zeta_{yx}} & & {V} \zeta_y  
\ar[dd]^{I\otimes u_y} \ar[rr]^{I\otimes {\zeta_{zy}}} & & 
V^2 \zeta_z \ar[dd]^{I^2\otimes u_z} 
\ar[rr]^{\FT \otimes I_{\zeta_z}} & & \ar[dd]^{I\otimes u_z} V \zeta_z \\
& & & & & & \\
\omega_x \ar[rr]_{\omega_{yx}} & & {V} \omega_y 
\ar[rr]_{I \otimes \omega_{zy}} & & V^2 \omega_z 
\ar[rr]_{\FT \otimes I_{\omega_z}} & & V \omega_z \\
}
$$
\subsubsection*{}

\subsubsection*{}
$$
\pic{
\zeta_x \ar[dd]_{u_x} \ar@{<-}[rr]^{\zeta_{xy}} & & {V} \zeta_y  
\ar[dd]^{I\otimes u_y} \ar@{<-}[rr]^{I\otimes {\zeta_{yz}}} & & 
V^2 \zeta_z \ar[dd]^{I^2\otimes u_z} 
\ar@{<-}[rr]^{F \otimes I_{\zeta_z}} & & \ar[dd]^{I\otimes u_z} V \zeta_z \\
& & & & & & \\
\omega_x \ar@{<-}[rr]_{\omega_{xy}} & & {V} \omega_y 
\ar@{<-}[rr]_{I \otimes \omega_{yz}} & & V^2 \omega_z 
\ar@{<-}[rr]_{F \otimes I_{\omega_z}} & & V \omega_z \\
}
$$
\subsubsection*{}

\noi and we can define the action of $\F_{\s}$ and of $\F_{\sb}$ on the arrows by

\bea
\F_{\s} (u_x,u_y,u_z) &=& (u_x,u_z) \\
\F_{\sb} (v_x,v_z) &=& (v_x, I \otimes v_z , v_z) \\
\eea
These functors satisfy the relations :

\bea
\F_{\s} \F_{\sb} (\xi_x,\xi_{zx},\xi_z) &=& (\xi_x,2\xi_{zx},\xi_z) \\
\F_{\sb} \F_{\s} ({\zeta_x},{\zeta_{yx}},{\zeta_y},{\zeta_{zy}},{\zeta_z}) 
&=& (\zeta_x,\zeta_{yx}\circ(A\otimes \zeta_{zy}),{V} \zeta_z,I \otimes I_{\zeta_z},\zeta_z) \\
\eea
If $(\alpha,\beta)\in\P_1$ is a generic 2-edge, then $\F_{\alpha\beta}$ acts locally as above 
without modifying the other entries.

For $p=0,\cdots,n$, let $\G_p=(\g_p,\cdots,\g_n)$ be the partial 2-path made of the last 
$(n-p+1)$ entries of $\G$ and let

\bea
\F_{\G_p} = \F_{\g_p\g_{p+1}}\circ\cdots\circ\F_{\g_{n-1}\g_n} 
&:& \F_{\g_n} \lra \F_{\g_p} \\
\eea
be the functor which maps the sections of $\f$ over $\g_n$ to sections over $\g_p$.
For example, we can choose $\g_n=(ab)$. Let us apply $\F_{\G_p}$ to the section 
$\zeta^n \in \ob (\F_{(ab)})$ defined by

\bea
\zeta^n_a &=& \zeta^n_b \ = \ V \\ 
\zeta^n_{ba} &=& F \\
\zeta^n_{ab} &=& \FT \\
\eea

\subsubsection*{}
$$
\zeta^n = \bigg( 
\pic{
{V} \ar@/_1pc/[rr]_{F} & & {V}^2 \ar@/_1pc/[ll]_{\FT} \\
}
\bigg)
$$
\subsubsection*{}

\subsubsection*{}
${\bf Theorem \  3}$ : \emph{$\F_\G$ multiplies the arrows of $\zeta^n$ by $K_T$ : }

\bea
\boxed{\F_{\G}(\zeta^n) = 
\bigg( 
\pic{
{V} \ar@/_1pc/[rr]_{K_T \, F} & & V^2 \ar@/_1pc/[ll]_{K_T \, \FT} \\
}
\bigg) \in \ob (\F_{(ab)}) } \\
\eea
\subsubsection*{}

\pr The inverse transport operator of $\zeta^p:=\F_{\G_p}(\zeta^n) \in \ob (\F_{\g_p})$ is

\bea
{\overline{T}}_{\g_p} (\zeta^p) = 
{\overline{T}}_{\g_p} \big( \F_{\g_p\g_{p+1}}\circ\cdots\circ 
\F_{\g_{n-1}\g_n} (\zeta^n) \big) : V^{\ell_p+1} \lra V  \\
\eea
By a decreasing induction on $p$, we have :

\bea
{\overline{T}}_{\g_p} (\zeta^p) 
= \zeta^n_{ab} \circ (\f_{\G_p}\otimes I) 
= \FT \circ (\f_{\G_p}\otimes I) \\
\eea

\subsubsection*{}
$$
\pic{
V^{\ell_p+1} \ar@/_1pc/[rrrr]_{{\overline{T}}_{\g_p} (\zeta^p) }
\ar@/^1pc/[rr]^{\f_{\G_p}\otimes I} & & {V}^2 \ar@/^1pc/[rr]^{\zeta^n_{ab}} & & V \\
}
$$
\subsubsection*{}

For $p=0$ :

\bea
\zeta^0_{ab} = {\overline{T}}_{\g_0} (\zeta^0) 
= \FT \circ (\f_{\G}\otimes I) = K_T\, \FT \\
\eea
Similarly, by using the direct transport operator, we obtain

\subsubsection*{}
$$
\pic{
V \ar@/_1pc/[rrrr]_{T_{\g_p} (\zeta^p) }\ar@/^1pc/[rr]^{\zeta^n_{ba}} 
& & V^2 \ar@/^1pc/[rr]^{\f_{{\widetilde{\G}}_p}\otimes I} & & V^{\ell_p+1} \\
}
$$
\subsubsection*{}

\bea
\zeta^0_{ba}=T_{\g_0} (\zeta^0)= ({\f}_{\widetilde{\G}}\otimes I) \circ F = K_T\, F \\
\eea

\hfill\qed

\subsubsection*{}

Once $\G=(\g_0,\cdots,\g_n)$ has been chosen, the sweeping process constructs 
a $V$-module, $\zeta_x={V}^{n_x}$, for each $x\in T_0$, and a morphism, $\zeta_{xy}$, 
for each oriented edge of $T$. Each integer $n_x$ depends only on the partial 2-path $\G_p$ 
which reaches $x$ first and not on the paths $\g_q$ with $q<p$. Similarly for each arrow, 
$\zeta_{xy}$. Therefore, we obtain a global section, $\zeta$, of $\f$ over $T$. 
More precisely, if we lift $T$ to a triangulation $\TT$ of the disk $D^2$ 
such that ${\TT}_{\vert\p D^2}$ be a pair of arcs both projected onto the base edge $(ab)$, 
then $\zeta$ is a global section of the pull-back of $\f$ to $\TT$.

If $\zeta_{xy}=0$ for some edge $(xy)$, then the transport operator along a path
$\g_p$ containing $(xy)$ vanishes, as well as the subsequent transport operators
and, at the end, we obtain $K_T=0$. Conversely, if $K_T=0$, then there exists an edge
(at least the last one) where $\zeta$ vanishes. Consequently, we have obtained
a geometric interpretation of the four-colour theorem in terms of sections of $\f$ :

\subsubsection*{}
${\bf Theorem \  4}$ : \emph{4CT $\iff ( \zeta_{xy}\neq 0 \quad \forall\, (xy) )$.}
\subsubsection*{}

{\bf Example : } Let us construct $\zeta$ on the octahedron.
Starting from $\zeta_a = \zeta_b = V $ and $\zeta^n_{ab} = \FT$, we have :

\bea
\zeta_e &=& V^2 \\
\zeta_{ae} &=& \FT \circ (\FT\otimes I) \\
\zeta_{eb} &=& I^2 \\
\zeta_d &=& V^3 \\
\zeta_{ad} &=& \FT \circ (\FT\otimes I) \circ (\FT\otimes I^2) \\
\zeta_{de} &=& I^3 \\
\zeta_f &=& V^2 \\
\zeta_{ef} &=& \FT \otimes I \\
\zeta_{fb} &=& I^2 \\
\zeta_{df} &=& (I\otimes \FT\otimes I) \circ (F\otimes I^2) \\
\zeta_c &=& V^3 \\
\zeta_{dc} &=& (I\otimes \FT\otimes I) \circ (F\otimes I^2)\circ (\FT\otimes I^3) \\
\zeta_{cf} &=& I^3 \\
\zeta_{ac} &=& \FT \circ (\FT\otimes I) \circ (\FT\otimes I^2)\circ (I^2\otimes \FT\otimes I) \\
& & \quad \circ (I\otimes F\otimes I^2)\circ (I\otimes\FT\otimes I^3)\circ (F\otimes I^3) \\
\zeta_{cb} &=& F\otimes I \\
\zeta_{ab}^0 &=& \FT \circ (\FT\otimes I) \circ (\FT\otimes I^2)
\circ (I^2\otimes \FT\otimes I)\circ (I\otimes F\otimes I^2) \\
& & \quad \circ (I\otimes\FT\otimes I^3)\circ (F\otimes I^3) \circ 
(I\otimes F\otimes I)\circ (F\otimes I) \\
&=& \FT \circ (\f_{\G}\otimes I) \\
&=& K_T\, \FT \\
\eea

\section{Conclusion and perspectives}

The classical approaches to the four-colour problem study the local form of a planar map 
to prove its global colourability. This suggests the existence of a cohomological interpretation 
of this property. In the present work, we have constructed a global section of a fibered category 
modeled on $\RG$ and proved that the validity of the four-colour theorem is equivalent to the fact 
that this section does not vanish. We hope that the present approach will be a first step toward 
an algebraic proof and the understanding of the four-colour theorem.

\vs

\end{document}